\documentclass{article}
\usepackage{lineno,hyperref}
\usepackage[american]{babel}
\usepackage{amssymb, amsthm,amsmath}
\usepackage{bm}
\usepackage{graphicx}
\usepackage[super]{cite}

\newtheorem{theorem}{Theorem}

\newtheorem{remark}[theorem]{Remark}

\usepackage{mathptmx}       
\usepackage{helvet}         
\usepackage{courier}        
\usepackage{type1cm}        
\usepackage{makeidx}         
\usepackage{graphicx}        
\usepackage{multicol}        
\usepackage[bottom]{footmisc}
\usepackage{amsmath}
\usepackage{amssymb}
\usepackage{amstext} 
\usepackage{float}
\allowdisplaybreaks
\makeindex             

\begin{document}
	\title{Multipatch ZIKV Model and Simulations}
	
	\author{A.~Sherly and W.~Bock\\[.3cm]
		{\small Technische Universit\"at Kaiserslautern,}\\{\small  Fachbereich Mathematik,}\\ {\small Gottlieb-Daimler-Stra{\ss}e 48,}\\ {\small 67663 Kaiserslautern, Germany}\\
		{\small E-Mail: $\{sherly,bock\}$@mathematik.uni-kl.de}\\[.2cm]
		}

	\maketitle


\begin{abstract}In this article we compare two multi-patch models for the spread of Zika virus based on an SIRUV model. When the commuting between patches is ceased we expect that all the patches follow the dynamics of the single patch model. We show in an example that the effective population size should be used rather than the population size of the respective patch.
\end{abstract}
\section{Introduction}\index{Multi-patch Zika Virus model}
Zika Virus belongs to the family \textit{Flaviviridae}, genus \textit{Flavivirus}. ZIKV disease is primarily vector-borne, which is transmitted by \textit{Aedes} mosquitoes\cite{Aedes_ZIKV}. This disease is also found to be sexually transmissible\cite{Sexually_transmitted}. 
Eventhough most patients show mild symptoms recent studies show that this virus attack results in neurological disorders like Guillain-Barré syndrome (GBS)\cite{Guillian_Barre}. Another important characteristic of this virus is its pathogenicity to fetuses causing  Microcephaly in newborn babies\cite{Microcephaly}.\\
The history of ZIKV disease known so far starts with the isolation of Zika virus from a rhesus monkey in Uganda around April 1947. There onwards it has spread across the world with the largest outbreak recorded in 2015-16 across South America\cite{ZIKV_SA,Aedes_ZIKV}. With no vaccines or medications found so far the disease spread can only be controlled by non-pharmaceutical interventions. Also increased international travel, evolution and mutation of viruses and their transmitting agents like mosquitoes, suitable environmental conditions etc lead to an increase in further outbreaks even in lesser probable places. The influence of human mobility plays an important role in transmitting diseases across continents. With more flight connectivity and affordable modes of transport disease transmission can also be faster. The primary objective of this study is to include spatial dependence to the mechanistic model of ZIKV spread. This is very relevant as the parameters involved in the model will be different for different places. So the dynamics will be exhibiting variations spatially. In this article we use an SIRUV model to describe the disease dynamics. This model divides the population into various compartments namely susceptible, infected and recovered. The interaction between various host and vector compartments, spread across different patches, is modeled using a coupling matrix and certain parameters.\\   
 We have discussed two models in section 2 and 3. The results of numerical simulations are provided in section 4. Comparing the two models exemplarily shows that the incorporation of the effective population size is crucial. While in a model, which just takes into account, the total population size of the patches, a decoupling does not lead to  the single patch dynamics, where as a model which incorporates the effective population size shows this desired property. 
\section{Multi-patch ZIKV model}\label{section2}\index{Multi-patch Zika Virus model}
 In this section we give a multi patch model for studing the ZIKV disease spread. Let the space domain be divided into small areas which we name as patches. The ZIKV model in a specific patch is also developed using different compartments. Here the host and vector population consists respectively of susceptible and infected compartments in each patch and we consider the recovered ones only in host population of each patch. We use either a subscript or a superscript ($i$, $j$ or $k$) to distinguish these compartments and the parameters patchwise. Let us first assume that the whole population is commuting between the patches and the rate of transition from patch $\left(i\right)$ to $
\left(j\right)$ be $p_{ij}$.
\begin{remark}
	The matrix $P$ with entries $p_{ij}$ is the residence time budgeting matrix. Here $p_{ij}$ represents the time spent by people in patch $i$ on average in patch $j$ in unit time\cite{P_matrix}. For example on average if a person in patch $i$ spent 8 hours in patch $j$, then $p_{ij} = \frac{8}{24}$, provided that unit time is one day.
\end{remark}
 We have deduced the following model from similar models in the literature used for other epidemiological studies\cite{Yashika}.
\begin{align*}
\begin{split}
&\frac{d S_{i}}{dt}=\mu_{i}\left( 1-S_{i}\right)-S_{i}\Bigg(\sum_{1\leq j\leq n}\beta_{vh}^{j}p_{ij}V_{j}+\Bigg(\sum_{1\leq j\leq n}\beta_{hh}^{i}\left(p_{ij}+p_{ji}\right){I_{j}}\Bigg.\Bigg.\Bigg.\Bigg.-\beta_{hh}^{i}p_{ii}{I_{i}}\Bigg)\Bigg)\\
&\frac{dI_{i}}{dt}=S_{i}\Bigg(\sum_{1\leq j\leq n}\beta_{vh}^{j}p_{ij}V_{j}+\Bigg(\sum_{1\leq j\leq n}\beta_{hh}^{i}\left(p_{ij}+p_{ji}\right){I_{j}}\Bigg.\Bigg.\Bigg.\Bigg.-\beta_{hh}^{i}p_{ii}{I_{i}}\Bigg)\Bigg)-\mu_{i}I_{i}-\gamma_{i}I_{i}\\
&\frac{d R_{i}}{dt}= \gamma_{i}I_{i}-\mu_{i}R_{i}\\
&\frac{d U_{i}}{dt}= \nu_{i}\left( 1-U_{i}\right) -\vartheta_{i}U_{i}\sum_{1\leq j\leq n}\frac{I_{j}}p_{ji}\\
&\frac{dV_{i}}{dt}=\vartheta_{i}U_{i}\sum_{1\leq j\leq n}{I_{j}}p_{ji}-\nu_{i}V_{i}.
\end{split}
\end{align*}
\section{Redefining the model for ZIKV}\label{section3}
Following some insights from \cite{Bichara} and \cite{Hethcote} we have developed a new model to describe the ZIKV disease spread. 
 In \cite{Hethcote}  a term called contact rate is clearly defined, which is the average number of adequate contacts per day of an infective person from patch $j$ with any individuals in patch $i$. 
With this in consideration we redefine the parameters used as follows\\
\\
$\alpha_j$ = number of infectious contacts that is happening per infected mostiquito per unit time with the people present in patch $j$.\\
$\beta_j$= number of infectious contacts that is happening per infective individual per unit time with the people present in patch $j$.\\
$\gamma_j$ = number of recoveries that is happening per unit time in patch $j$.\\
$\vartheta_j$ = number of infective contacts that is happening per infected human with mosquitoes in patch $j$ in unit time.
\\\\
Let us focus on patch $j$ and see how many susceptibles from patch $i$ is infected in patch $j$. By the definition of $\alpha_j$ the number of people getting into adequate contacts with the mosquitoes in patch $j$ is given by $\alpha_j \mathcal{V}_j$. Now the total number of people who were present in patch $j$ is given by $\sum_{k=1}^n p_{kj}N_k$. Let us call this the effective population in patch $j$. Also the effective population of susceptibles in patch $j$ is $\sum_{k=1}^{n} p_{kj}\mathcal{S}_k$. 
Among which $p_{ij}\mathcal{S}_i$ are coming from patch $i$. The number of susceptibles from patch $i$ who get infected in patch $j$ due to mosquitoes is given by $$\alpha_j\mathcal{V}_j\frac{p_{ij}\mathcal{S}_i}{\sum_{k=1}^n p_{kj}N_k}.$$ Now we focus on the infections between humans. The number of infections happening in patch $j$ in unit time due to human-human interactions is given by $\beta_jI_{eff}$,  where $I_{eff}$ is the effective number of infected people who came to patch $j$ in unit time which is given by $\sum_{k=1}^{n} p_{kj}\mathcal{I}_k$. The total number of infections happening in patch $j$ is given by $\beta_j\sum_{k=1}^{n} p_{kj}\mathcal{I}_k$ out of which the number of infections happened to the susceptible people of patch $i$ is $$\beta_j\sum_{k=1}^{n} p_{kj}I_k\frac{p_{ij}\mathcal{S}_i}{\sum_{k=1}^n p_{kj}N_k}.$$
Now we have to introduce fractions by normalising each compartmental values. 
\begin{center}
	\begin{tabular}{c c c c c}
		${S}_i$&${I}_i$ &${R}_i$&${U}_i$ &${V}_i$\\\hline\noalign{\smallskip} $\frac{\mathcal{S}_i}{N_i}$&$\frac{\mathcal{I}_i}{N_i}$&$\frac{\mathcal{R}_i}{N_i}$
		&$\frac{\mathcal{U}_i}{M_i}$&$\frac{\mathcal{V}_i}{M_i}$
	\end{tabular}
\end{center}
The following system of ODEs describe disease spread in each patch $i$
\begin{align*}
\frac{d{S}_{i}}{dt}&= \mu_i(1-{S}_i)-{\sum_{j=1}^{n}\alpha_jM_j{V}_j\frac{p_{ij}{S}_i}{\sum_{k=1}^n p_{kj}N_k}}-{\sum_{j=1}^{n}\beta_j\sum_{k=1}^{n} p_{kj}N_k{I}_k\frac{p_{ij}{S}_i}{\sum_{k=1}^n p_{kj}N_k}}\\
\frac{d{I}_{i}}{dt}&= -(\gamma_i+\mu_i){I}_i+{\sum_{j=1}^{n}\alpha_jM_j{V}_j\frac{p_{ij}{S}_i}{\sum_{k=1}^n p_{kj}N_k}}+{\sum_{j=1}^{n}\beta_j\sum_{k=1}^{n} p_{kj}N_k{I}_k\frac{p_{ij}{S}_i}{\sum_{k=1}^n p_{kj}N_k}}\\
\frac{d{R}_{i}}{dt}&= \gamma_i{I}_i-\mu_i{R}_i\\
\frac{d{U}_{i}}{dt}&= \nu_i(1-{U}_i)-\vartheta_i\frac{{U}_i}{M_i}{\sum_{k=1}^n p_{ki}N_k{I}_k}\\
\frac{d{V}_{i}}{dt}&= -\nu_i {V}_i+\vartheta_i\frac{{U}_i}{M_i}{\sum_{k=1}^n p_{ki}N_k{I}_k}.
\end{align*}
\section{Comparison of both models in three-patch scenario}
In a case where $n=3$ we numerically simulated both the models and  compared the results. We obtained the influence of the residence time budgeting matrix on the multi-patch model. Here we restrict ourselves to consider three patches with the same set of parameters and population sizes. The movements between these three patches are defined using the residence time budgeting matrix $P$. The question is how far does the dynamics deviate from the single patch case, when the movement between the patches is controlled using the $p_{ij}$ values. We use the parameters 
and population sizes, as given in table \ref{table1}, for the numerical simulation. We are studying two cases- the three patches being coupled and completely decoupled respectively. For the first case 
\begin{align}\label{P_matrix}
P =
 \begin{bmatrix}
 0.2 &0.7&0.1\\
 0.5&0.1&0.4\\
 0.3&0.6&0.1
  \end{bmatrix}
  \end{align}
\begin{table} 
\begin{center}
	\caption{Parameters and Population Sizes}
	{\footnotesize
		\begin{tabular}{p{2cm}p{1cm}p{1cm}p{1cm}p{1cm}p{1.5cm}p{1.5cm}}
			\hline\noalign{\smallskip}
			$\mu$ & $\alpha$ &$\beta$&$\vartheta$&$\nu$&N&M \\
\hline
\\
			10/(1000*365)          
			&0.008						
			& 0.01                    			
			&0.4                    			
			& 1/14                  			
			& 20000
			& 100000\\
			\noalign{\smallskip}\hline\noalign{\smallskip}
	\end{tabular}}
\label{table1}
\end{center}
\end{table}
\begin{figure}[ht]
	\centering
	\includegraphics[width=0.7\linewidth]{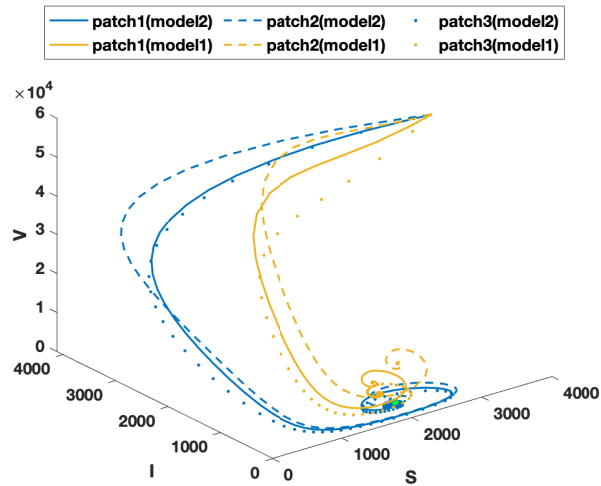}
	\caption{Phase portrait for three patches using model 1(Section \ref{section2})  and model 2(Section \ref{section3}) for the case where the patches are coupled using the matrix $P$ from \eqref{P_matrix}.}
	\label{fig:fig2phaseportrait}
\end{figure}
\begin{figure}[H]
	\centering
	\includegraphics[width=0.7\linewidth]{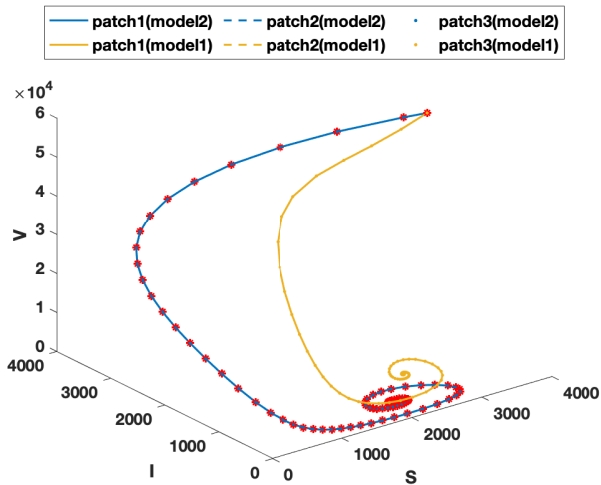}
	\caption{For the same set of parameters as in Figure \ref{fig:fig2phaseportrait} when $P$ is set to identity matrix we see the given results where the red starred curve is the phase portrait of the single patch model}
	\label{fig:fig2pzero}
\end{figure}
The dynamics was supposed to be similar for the single patch and multi-patch models for the case $P=I$. But we have not seen this property for the old model.
\section{Conclusion}
In this study we have considered two different models to describe the dynamics of ZIKV spread. We compared the two models to identify the suitable model. When the commuting between patches is ceased we expect that all the three patches follow the dynamics of the single patch model. The first model failed to satisfy this condition where as the second model was successfully exhibiting this property. This gives rise to a more thorough study of the second model in a forthcoming work.


\begin{thebibliography}{1}
	
	\bibitem{Aedes_ZIKV}
	E.~B. Kauffman and L.~D. Kramer, ``{Zika Virus Mosquito Vectors: Competence,
		Biology, and Vector Control},'' {\em The Journal of Infectious Diseases},
	vol.~216, pp.~S976--S990, 12 2017.
	
	\bibitem{Sexually_transmitted}
	P.~S. Mead, S.~L. Hills, and J.~T. Brooks, ``Zika virus as a sexually
	transmitted pathogen,'' {\em Current Opinion in Infectious Diseases},
	vol.~31, no.~1, 2018.
	
	\bibitem{Guillian_Barre}
	L.~Barbi, A.~V.~C. Coelho, L.~C. A.~d. Alencar, and S.~Crovella, ``Prevalence
	of guillain-barr{\'e}syndrome among zika virus infected cases: a systematic
	review and meta-analysis,'' {\em The Brazilian Journal of Infectious
		Diseases}, vol.~22, no.~2, pp.~137--141, 2018.
	
	\bibitem{Microcephaly}
	A.~Q.~C. Araujo, M.~T.~T. Silva, and A.~P. Q.~C. Araujo, ``{Zika
		virus-associated neurological disorders: a review},'' {\em Brain}, vol.~139,
	pp.~2122--2130, 06 2016.
	
	\bibitem{ZIKV_SA}
	A.~S. Fauci and D.~M. Morens, ``Zika virus in the americas ---yet another
	arbovirus threat,'' {\em New England Journal of Medicine}, vol.~374,
	pp.~601--604, 2021/08/27 2016.
	
	\bibitem{P_matrix}
	P.~Heidrich, Y.~Jayathunga, W.~Bock, and T.~G{\"o}tz, ``Prediction of dengue
	cases based on human mobility and seasonality---an example for the city of
	jakarta,'' {\em Mathematical Methods in the Applied Sciences}, vol.~n/a,
	2021/08/27 2021.
	
	\bibitem{Yashika}
	W.~Bock and Y.~Jayathunga, ``Optimal control and basic reproduction numbers for
	a compartmental spatial multipatch dengue model,'' {\em Mathematical Methods
		in the Applied Sciences}, vol.~41, pp.~3231--3245, 2021/08/27 2018.
	
	\bibitem{Bichara}
	D.~Bichara and A.~Iggidr, ``Multi-patch and multi-group epidemic models: a new
	framework,'' {\em Journal of Mathematical Biology}, vol.~77, no.~1,
	pp.~107--134, 2018.
	
	\bibitem{Hethcote}
	H.~W. Hethcote, ``The mathematics of infectious diseases,'' {\em SIAM Review},
	vol.~42, no.~4, pp.~599--653, 2000.
	
\end{thebibliography}
\end{document}